\documentclass[12pt]{article}
\usepackage{amssymb}
\usepackage{epsfig}
\usepackage{amsmath}
\usepackage{amsthm}
\usepackage{subfigure}
\usepackage{url}
\usepackage{color}
\usepackage{graphicx}
\DeclareMathOperator{\transverse}{\cap\kern-7.75pt\top}
\newtheorem{con}{Conjecture}[section]

\newcommand{\argmin}{\operatorname{argmin}}
\newcommand{\spt}{\operatorname{spt}}

\newcommand{\Per}{\operatorname{Per}}
\newtheorem{thm}[con]{Theorem}
\newtheorem{lem}[con]{Lemma}
\newtheorem{deff}[con]{Definition}

\newtheorem{prop}[con]{Proposition}
\newtheorem{rem}[con]{Remark}

\newcommand{\myell}{\,\hbox{\vrule height 6pt depth 0pt
\vrule height 0.4pt depth 0pt width 5pt}\,}

\usepackage[bookmarks=true, % true means bookmarks
            bookmarksnumbered=false, % left window are numbered
            bookmarksopen=false, % true means only level 1 are displayed. 
            colorlinks=true, 
            linkcolor=red]{hyperref}

\begin{document}

\title{{\bf $L^1$TV computes the flat norm
    for boundaries}} 

\author{Simon P. Morgan$^*$ and Kevin R. Vixie$^{**}$}

\date{$^*$University of Minnesota and $^{**}$Los Alamos National Laboratory}

\maketitle

\begin{abstract}
  We show that the recently introduced $L^1$TV functional can be used
  to explicitly compute the flat norm for co-dimension one boundaries.
  Furthermore, using $L^1$TV, we also obtain the flat norm
  decomposition.  Conversely, using the flat norm as the precise
  generalization of $L^1$TV functional, we obtain a method for
  denoising non-boundary or higher co-dimension sets.  The flat norm
  decomposition of differences can made to depend on scale using the
  \emph{flat norm with scale} which we define in direct analogy to the
  $L^1$TV functional. We illustrate the results and implications with
  examples and figures.
\end{abstract}

\section{Introduction}
\label{sec:intro}
In this research announcement, we point out that the
$L^1$TV functional, introduced and studied in the continuous setting
in \cite{chan-2005-1, vixie-2006-1, allard-2006-1, yin-2006-1} and
earlier in the discrete setting in
\cite{alliney-1997-1,nikolova-2003-1}, provides a convenient way of
computing both the value of, and the optimal decomposition required
by, the \emph{flat norm} from geometric measure theory.

The $L^1$TV functional was introduced as an improvement to the now
classic Rudin-Osher-Fatemi total variation based denoising~\cite{rudin-1992-1}. Among its
many nice properties are the way it handles binary images, making it
useful for shape processing. Theoretically speaking, the clean
geometric structure of the functional and its minimizers is very
attractive. See \cite{vixie-2006-1, allard-2006-1} for details.

Joan Glaun\`{e}s was, as far as we know, the first to suggest and use
the \emph{flat norm} from geometric measure theory as a distance in
shape space. In his dissertation~\cite{glaunes-2005-1} and a couple of
application papers~\cite{vaillant-2005-1,glaunes-2007-1} with
collaborators, the dual formulation of the flat norm is used to
compute distances between shapes.

{\bf Acknowledgements}: It is a pleasure to acknowledge Bill Allard
and Bob Hardt for useful discussions and Joan Glaun\`{e}s and Sarang
Joshi for inspiring this work.  Additionally, the second author
acknowledges useful discussions with Selim Esedoglu.

\section{$L^1$TV gives the Flat Norm}
\label{sec:l1tv-flat}
The $L^1$TV functional introduced and studied in~\cite{chan-2005-1} is given by:
\begin{equation}
  \label{eq:l1tv}
  F(u) = \int_{\Bbb{R}^n} |\nabla u| dx + \lambda \int_{\Bbb{R}^n} |u-u_0| dx
\end{equation}
where $u$ and $u_0$ are functions from $\Bbb{R}^n$ to $\Bbb{R}$ with
$u_0$ being the input function.  In image analysis applications, $n$
is usually 2 and $u_0$ is the measured image intensity function. The
optimal $u$ minimizing~(\ref{eq:l1tv}) can be thought of as a
\emph{denoised} version of $u_0$. Typically, one chooses the parameter
$\lambda$ based on noise levels since this choice selects the scale
below which features or oscillations are ignored. In the case that the
input function is binary, Chan and Esedoglu observe that the
functional reduces to:
\begin{equation}
  \label{eq:set-l1tv}
  F_{CE}^{\lambda}(\Sigma) = \Per(\Sigma) + \lambda |\Sigma\vartriangle\Omega|.
\end{equation}
where the binary minimizer is $\chi_{\Sigma}$, the binary input is
$\chi_{\Omega}$, $\Per(\Sigma)$ is the perimeter of $\Sigma$, and
$\vartriangle$ denotes the symmetric difference. Of course $\chi_{E}$
denotes the characteristic function on $E$, with a value of 1 on $E$
and 0 on the complement of $E$. Now, let $\Sigma(\Omega,\lambda)$ be a
binary minimizer of ~(\ref{eq:set-l1tv}). I.e.
\begin{equation}
  \label{eq:minim-l1tv}
  \Sigma(\Omega,\lambda) \equiv \argmin F_{CE}^{\lambda}(\Sigma) = \Per(\Sigma) + \lambda |\Sigma\vartriangle\Omega|.
\end{equation}
For our convenience, we record the optimal decomposition of $\Omega$
into $\{\Sigma(\Omega,\lambda)$ and
$(\Sigma(\Omega,\lambda)\vartriangle\Omega)\}$ as the pair $\{\partial
\Omega, \Sigma(\Omega,\lambda)\vartriangle\Omega\}$.

In what follows, we use the notions of \emph{current}, \emph{mass} and
\emph{$\partial$} and supporting ideas from geometric measure theory.
We introduce and explain these in some detail in the Appendix.
Informally, one can gain much by thinking of the \emph{n-current}
$T_E$ as an orientable n-submanifold or n-rectifiable set $E \subset
\Bbb{R}^{n+k}$ with an orientation, of the \emph{mass} $M(T_E)$ as the
$n$-dimensional volume of $E$, and of $\partial T$, the boundary of
$T$, as the oriented boundary of $E$ with orientation imposed by the
orientation of $E$. We sometimes omit the subscript indicating the
support of the current, referring to currents $T$ and $S$.  For those
without experience with currents, we suggest focusing the examples in
the the appendix.

The \emph{flat norm} of an $n$-current $T$, denoted by $\Bbb{F}(T)$,
is given by
\begin{equation}
  \label{eq:flat-def}
\Bbb{F}(T) \equiv \min_{S} \{\mathbf{M}(S) + \mathbf{M}(T - \partial S)\}) 
\end{equation}
where $S$ varies over $n+1$-currents and $\mathbf{M}$ is the
\emph{mass} of the indicated currents. We refer to $\{T,S\}$ as the
flat norm induced, optimal decomposition. Now for results. 

\begin{thm}
    For the current $T_{\partial \Omega}$, 
\begin{equation}
  \label{eq:flat1}
  \Bbb{F}(T_{\partial \Omega}) = F_{CE}^{1}(\Sigma(\Omega,1))
\end{equation}
and
\begin{equation}
  \label{eq:flat2}
 \{T_{\partial \Omega}, S_{\Sigma(\Omega,1)\vartriangle\Omega}\}
\end{equation}
is the optimal decomposition required by the flat norm. 
\label{thm:litvfn} 
\end{thm}
\noindent Theorem~\ref{thm:litvfn} says that L1TV computes the flat
norm. This relation between the flat norm and the L1TV functional
immediately suggests a very useful generalization of the flat norm.
\begin{deff}[{\bf Flat Norm With Scale}]
  \begin{equation}
    \Bbb{F}_{\lambda}(T) \equiv \min_{S} \{\lambda\mathbf{M}(S) + 
            \mathbf{M}(T - \partial S)\})    
  \end{equation}
\end{deff}
\noindent Theorem~\ref{thm:litvfn} is then simply a special case of
the next theorem.
\begin{thm}
    For the current $T_{\partial \Omega}$, 
\begin{equation}
  \label{eq:flat1-param}
  \Bbb{F}_{\lambda}(T_{\partial \Omega}) = F_{CE}^{\lambda}(\Sigma(\Omega,\lambda))
\end{equation}
and
\begin{equation}
  \label{eq:flat2-param}
  \{T_{\partial \Omega}, S_{\Sigma(\Omega,\lambda)\vartriangle\Omega}\}
\end{equation}
is the optimal decomposition that the flat norm with scale requires. 
\label{thm:litvfn-ws}  
\end{thm}

\begin{proof}[{\bf Proof of Theorem~\ref{thm:litvfn-ws}}]
   The proof is really simply checking
that the picture one can draw holds after the definitions of
\emph{mass} ($\mathbf{M}$) and \emph{perimeter} ($\Per$) are used to
translate the picture into analytic terms. Very briefly, we have:
\begin{eqnarray*}
  \label{eq:proof1}
  F_{CE}^{\lambda} (\Sigma) & = & \lambda |\Sigma\vartriangle\Omega| 
                              + \Per(\Sigma)\\
                           & = & \lambda M(S_{\Sigma\vartriangle\Omega}) 
                              + M(T_{\partial \Sigma})\\
                           & = & \lambda M(S_{\Sigma\vartriangle\Omega}) 
                                + M(T_{\partial \Omega} 
                                - \partial S_{\Sigma\vartriangle\Omega})
\end{eqnarray*}
Figure~\ref{fig:symmetric-diff-picture} illustrates this pictorially.
\begin{figure}[htp!]
\begin{center}
\input{symmetric-diff-picture.pstex_t}
\end{center}
\caption{In this figure we illustrate the translation of the $L^1$TV
  view to the flat norm view.  The perimeter of $\Sigma$ becomes the
  mass of $T_{\partial\Omega}-\partial S_{\Sigma\vartriangle\Omega}$
  and the volume of $\Sigma\vartriangle\Omega$ becomes the mass of
  $S_{\Sigma\vartriangle\Omega}$.  Note: this figure does not depict a
  minimizer. Rather, we depict $\Omega$ and any candidate $\Sigma$.}
\label{fig:symmetric-diff-picture}
\end{figure}
\end{proof}

Our final observation is that the \emph{flat norm with scale} gives us
the same decomposition as we would get if we first scaled $T$,
computed the \emph{flat norm} decomposition and then reversed the
scaling. More precisely
\begin{lem}
  Denote the optimal $\Bbb{F}_{\lambda}$ decomposition by
  $\{T,S\}_{\lambda}$. Then
 \begin{equation}
   \label{eq:fn-scaling}
   \{T,S\}_{\lambda} = d_{\frac{1}{\lambda}\#}\{d_{\lambda\#}(T),d_{\lambda\#}(S)\}_1
 \end{equation}
 where $d_\lambda$ denotes the $\lambda$-dilation of $\Bbb{R}^m$, and
 $d_{\lambda\#}(T_M) = T_{d_{\lambda}(M)}$.
\end{lem}
\begin{proof}
The minimizing decomposition $\{T,S\}_{\lambda}$
which minimizes $\lambda\mathbf{M}(S) + \mathbf{M}(T - \partial S)$
also minimizes $\lambda^{n}\mathbf{M}(S) + \lambda^{n-1}\mathbf{M}(T -
\partial S) = \mathbf{M}(d_{\lambda\#}S) + \mathbf{M}(d_{\lambda\#}(T -
\partial S))$. Therefore, to get the minimizer for $\Bbb{F}_{\lambda}$
run the optimization required by $\Bbb{F}_{1}$ using $d_{\lambda\#}(T)$
as input, and then contract with $d_{\frac{1}{\lambda}\#}$.  
\end{proof}
We now discuss applications and examples with pictures which should
clarify things for those with less exposure to geometric measure
theory.

\section{Applications and Illustrating Examples}
\label{sec:examples}

The value of the above results is fully realized by exploring
their use in applications to images and shapes. As observed by
Glaun\`{e}s and others, the flat norm is a very natural candidate for
distances between shapes. We now very briefly explore and illustrate
applications of the above observations.

\subsection{Generalized Flat Norm: flat norm with scale}
\label{sec:gfn-fnws}

\begin{figure}[htp!]
\begin{center}
\input{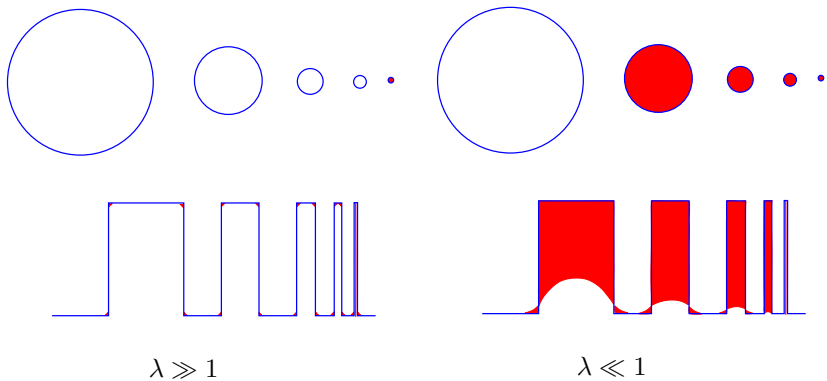}
\end{center}  
\caption{Flat norm with scale decompositions as a function of
  $\lambda$: The larger $\lambda$ is, the less smoothing there is. Let
  $T$ denote the blue input 1-currents and $S$ denote the red optimal
  2-currents. Then in fact, $\lambda$ is the bound on the
  curvature allowed in $T-\partial S$,
  see~\cite{allard-2006-1,vixie-2006-1} for details. Accordingly, the
  curves $T - \partial S$ are allowed greater curvature in the
  curves on the left than is permitted in those on the right.}
  \label{fig:scale-dependence-picture}
\end{figure}

As mentioned in section~\ref{sec:l1tv-flat},  by letting the
$\lambda \neq 1$, one has a natural way to
vary the intrinsic scale in the flat norm. 
\begin{eqnarray*}
  \label{eq:lambda-f}
  \Bbb{F}_{\lambda}(T_{\partial E}) & \equiv & \min_{S} \{\lambda \mathbf{M}(S) 
  + \mathbf{M}(T_{\partial E} - \partial S)\} \\
  &   =    & F_{CE}^{\lambda}(E) 
\end{eqnarray*}
where $T_{\partial E}$ represents $\partial E$. The
main point here is that this is easy to compute, given the connection
to the $L^1$TV functional. Varying $\lambda$ gives us the ability to
choose what scale is big and worth keeping. See
Figure~\ref{fig:scale-dependence-picture}.

\subsection{Flat norm via $L^1$TV}
\label{sec:fnv-l1tv}

\begin{figure}[htp!]
    \begin{center}
      \input{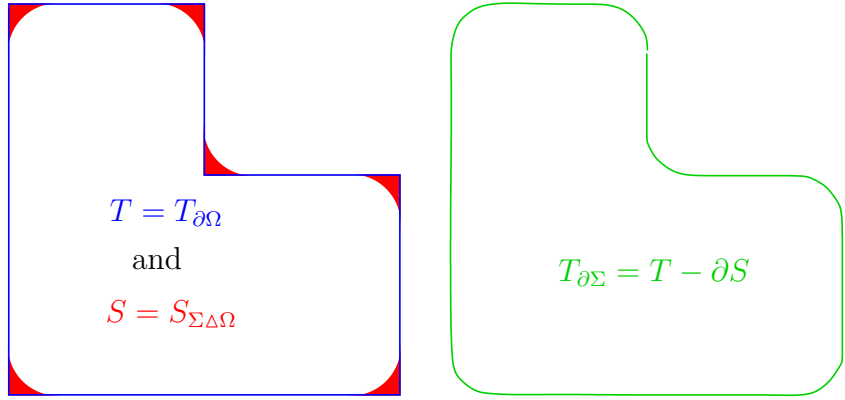}
    \end{center}
    \caption{A simple reinterpretation of the $L^1$TV input and
      minimizer gives us the flat norm of $T_{\partial\Omega}$ and the
      decomposition into the diminished boundary $T_{\partial\Sigma}$
      and the symmetric difference current
      $S_{\Sigma\vartriangle\Omega}$.}
    \label{fig:L1-through-flat-decomp}
  \end{figure}

We can use the $L^1$TV functional to very easily calculate \emph{both}
the flat norm of differences between surfaces which are boundaries
\emph{and} the optimal decomposition of that difference into surface
and area parts.  In Figure~\ref{fig:L1-through-flat-decomp}, the
decomposition of a boundary current $T_{\partial\Omega}$ into the
diminished boundary $T_{\partial\Sigma}$ and the symmetric difference
current $S_{\Sigma\vartriangle\Omega}$ is illustrated.

\subsection{$L^1$TV by the dual form of the flat norm}
\label{sec:l1tv-btdf}

The dual formulation of the flat norm can be used to compute $L^1$TV
minimizers. In what follows, we define $\spt T$ to be the support of
$T$.  The following results establish that maximizing forms, or in
some cases, maximizing sequences of forms, contain the decomposition
into $S$ and $T-\partial S$.

\begin{prop}
Suppose that $T(\phi )=F(T)=\underset{S}{min}(M(S)+M(T-\partial S))$,
where $\phi$ is a smooth, compactly supported $n$-form satisfying
$|\phi| \leq 1, |d\phi| \leq 1,$ and $T$ is an $n$-rectifiable current
with density $\Theta =1$ $\mathcal{H}^{n}$ almost everywhere on $\spt{T}$, then on
the support of $S,$ $\left\vert d\phi \right\vert =1,$ and on the
support of $T-\partial S$, $\left\vert \phi \right\vert =1.$
\label{prop:maxform}
\end{prop}

\begin{proof}
For a minimizing choice of S, we have that 
\begin{eqnarray*}
  T(\phi ) & = & F(T)=M(S)+M(T-\partial S)\\
  T(\phi ) & = & \partial S(\phi )+(T-\partial S)(\phi )
\end{eqnarray*}
and 
\begin{equation*}
\partial S(\phi) =  S(d\phi)
\end{equation*}
so that
\begin{equation}
S(d\phi) + (T-\partial S)(\phi ) =  M(S)+M(T-\partial S)
\end{equation}

We know $M(S) = \int 1 d\mathcal{H}^{n+1}\myell\spt{S}$ and $M(S) =
\int <d\phi,\vec{S}> d\mathcal{H}^{n+1}\myell\spt{S}$, $|\vec{S}| = 1$
and $|d\phi| \leq 1$ on $\spt{S}$. Similarly for $T - \partial S$.
We immediately get that $|\phi| = 1$, $\mathcal{H}^{n}
\myell\spt{(T - \partial S)} $ almost everywhere and $|d\phi| = 1$,
$\mathcal{H}^{n+1}\myell\spt{S}$ almost everywhere.
\end{proof}

\begin{figure}[htp!]
  \begin{center}
    \input{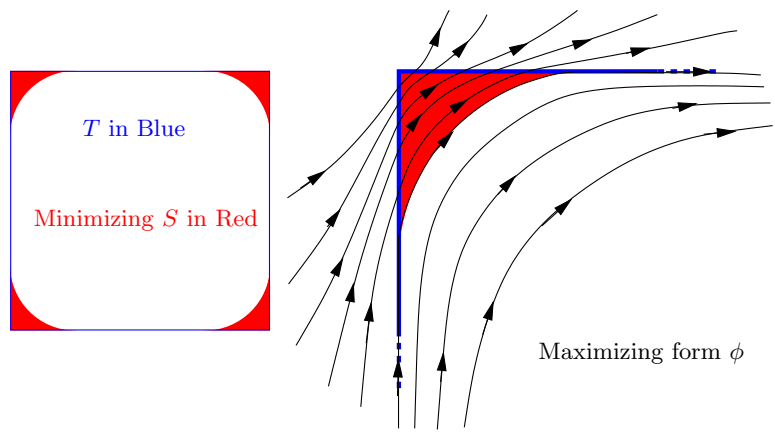}
  \end{center}
  \caption{Example of a maximizing form for a square $T$ as input. On
    the right hand side of the figure, details of the form are shown. We
    visualize the form as a vector field. On $T-\partial S$, $|\phi|
    = 1$ everywhere and on $S$ (red), $|d\phi|=1$ everywhere, while
    off of these sets, both $|\phi| < 1$ and $|d\phi| < 1$.}
  \label{fig:square-form-example}
\end{figure}

See Figure~\ref{fig:square-form-example} for an example maximizing
form. Note that when $S$ is not unique, this proposition implies that
$|d\phi| = 1$ on the union of supports of all possible minimizing $S$'s.
If we do not have a maximizing form, we have the following modified
proposition together with the fact that there will be sequence of
forms $\phi_i$ such that $T(\phi_i)\rightarrow_{i\rightarrow\infty} F(T)$.
This easily yields
\begin{prop} Suppose that $T(\phi_i)\rightarrow F(T) = M(S) +
  M(T-\partial S)$. Then there is a
  subsequence of $\phi_i$, $\phi_{i_k}$ such that
  $|\phi_{i_k}|\rightarrow 1$ $\mathcal{H}^{n} \myell\spt{(T-\partial S)} $ almost
  everywhere and $|d\phi|\rightarrow 1$
  $\mathcal{H}^{n+1}\myell\spt{S}$ almost everywhere.
\label{prop:maxseq}
\end{prop}

This modified proposition is necessary since there are easily
constructed examples having no smooth maximizing form. In fact, the
example shown in Figure~\ref{fig:square-form-example} is not actually
smooth. The optimizing $d\phi$ we show is actually Lipschitz, so it
can be arbitrarily well approximated by smooth forms even though it is
not itself smooth. The non-smoothness originates at the points of
$T-\partial S$ where the circular arcs join the sides of the square
tangentially. At these points, the boundary is merely $C^{1,1}$.

\begin{rem}
  Note that Propositions~\ref{prop:maxform} and~\ref{prop:maxseq} are
  true if $\phi$ and $d\phi$ are merely measurable wrt the measures
  $\mathcal{H}^{n} \myell\spt{(T-\partial S)}$ and
  $\mathcal{H}^{n+1}\myell\spt{S}$.
\end{rem}

We now state an outline of a reasonable set of steps that would be very
useful for computations using the dual formulation of the flat norm.

\begin{con} Suppose that $T$ is integer density rectifiable with
  rectifiable boundary $\partial T$ and that both $T$ and $\partial T$
  have finite mass.  Define $\Phi$ to be the collection of all
  Lipschitz forms $\phi$ maximizing $T(\phi)$ and satisfying $|\phi|
  \leq 1$ and $|d\phi| \leq 1$. Define $X$ to be the closure of the
  set on which $|d\phi|=1$ for every $\phi\in\Phi$. Finally, define
  $\Bbb{S}$ to be the collection of all optimizing currents $S$ such
  that $\Bbb{F}(T)= M(S) + M(T-\partial S)$. Then
\begin{enumerate}
\item[(a)] There is an $S \in \Bbb{S}$ such that both $S$ and
  $\partial S$ are integer density rectifiable.
\item[(b)] $\Phi \neq \emptyset$
\item[(c)] $X = \underset{S\in\Bbb{S}}{\bigcup} \;\;\spt(S)$
\item[(d)] $|\phi| = 1$ on $\spt(T-\partial S)$.
\end{enumerate}
\end{con} Part (a) says that although we are not constraining the
minimizing currents to be rectifiable, there is at least one in the
set of minimizers that is. Part (b) says that there are always
maximizing forms if we permit them to be merely Lipschitz instead of
smooth. Part (c) enables us to see where the possible locations for an
$S$ might be and finally, part (d) gives us $T-\partial S$.

\begin{rem}
  We do not expect the refinement and proof of the conjecture to be
  simple. The few previous works in this direction, for example
  Federer's 1974 paper~\cite{federer-1974-1} ``Real Flat Chains,
  Cochains and Variational Problems'', are rather technical in nature.
  Notice also that this conjecture is only necessary for a rigorous
  foundation to the use of the dual formulation in the computation of
  the flat norm decomposition. Direct optimization over rectifiable
  currents needs nothing from this conjecture for its justification.
\end{rem}

\begin{figure}[htp!]
  \centering
  \input{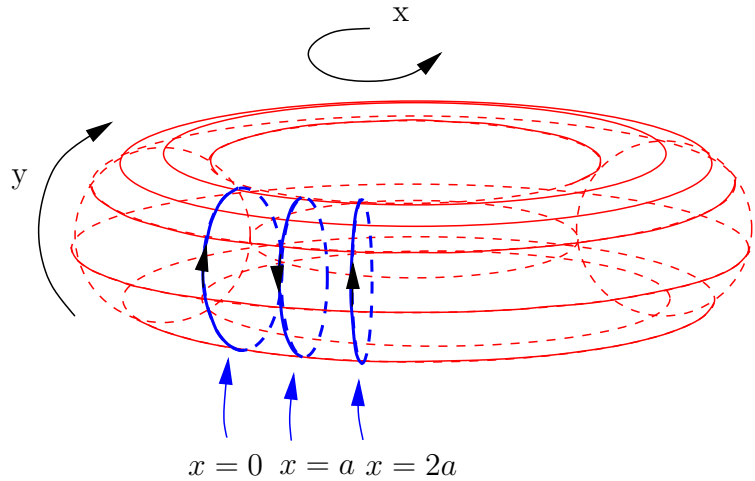}
  \caption{The 3 circles example. The metric is chosen so that the
    circles are parallel to each other. Optimal $S$ can be either the
    region between the $x=0$ and $x=a$ circles or the region between
    the $x=a$ and $x=2a$ circles.}
  \label{fig:3-lines-example}
\end{figure}

\begin{figure}[htp!]
  \centering
  \input{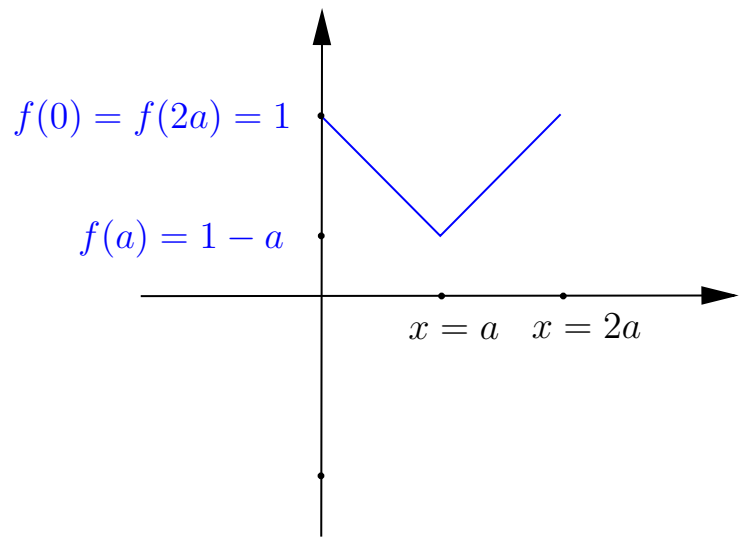}
  \caption{A maximizing form for the 3 circles example. We plot $f$
    for the form given by $f(x)dy$. We are forced to use $\alpha = 1$
    and the Lipschitz form plotted here.}
  \label{fig:3-lines-example-function}
\end{figure}

A very simple example where Lipschitz forms are necessary and
sufficient for optimality is the case in which the current is three
equally spaced circles on a torus. See Figure
~\ref{fig:3-lines-example}. Note that the metric on the torus is
chosen such that the circles are parallel.
Figure~\ref{fig:3-lines-example-function} shows the $f$ of a Lipschitz
maximizing form $f(x)dy$.  In the case shown of equal spacing
between circles, we can't maximize with a smooth form and a maximizing
sequence must approach the form plotted in
Figure~\ref{fig:3-lines-example-function}. In this case, the region
between the $x=0$ and $x=a$ circles or the region between the $x=a$
and $x=2a$ circles is the optimal $S$. This non-uniqueness is taken
into account in the above conjecture.
Finally, the above propositions, example and conjecture have obvious
analogs for $\Bbb{F}_{\lambda}$, the flat norm with scale.

\subsection{$L^1$TV for co-dimension $> 1$}
\label{sec:l1tv-cdgt1}

\begin{figure}[htp!]
    \begin{center}
      \input{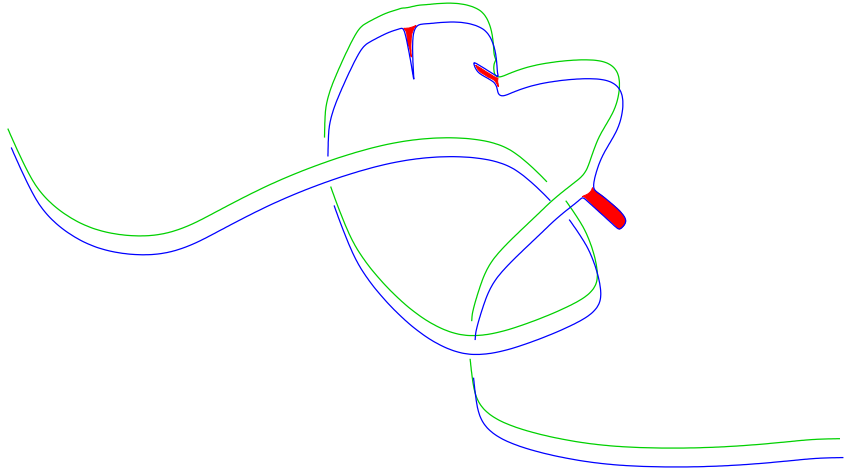}
    \end{center}
    \caption{The green curve is the denoised version of the blue, where
      we have translated the green to make visualization easier.}
    \label{fig:co-dimension-2}
  \end{figure}

Computing the flat norm decomposition for arbitrary currents permits
us to extend the $L^1$TV denoising to sets which are not boundaries or
have co-dimension greater than 1. One approach is to use the dual
formulation of the flat norm. This depends on extracting the optimal
decomposition from the optimizing form, as discussed in the previous
subsection. Another approach is to directly optimize over currents. We
are currently developing both approaches.
Figure~\ref{fig:co-dimension-2} schematically illustrates the
decomposition of a 1-current in 3D that results when the flat norm is
computed. This example actually illustrates both of the useful
generalizations possessed by the flat norm decomposition:
regularization or denoising of higher co-dimension and non-boundary
subsets. Notice that the use of the flat norm with scale permits us to
choose what scale is small and therefore greatly diminished, and what
scales are large and therefore preserved.  In the case of sets which
are co-dimension 1 boundaries, we know that in a very precise sense,
the regularized surface given by $\spt(T-\partial S)$ is the best
$\lambda$-curvature approximation to $T$.
See~\cite{vixie-2006-1,allard-2006-1} for details.

\subsection{Shape Statistics}
\label{sec:shape-statistics}

\begin{figure}[htp!]
  \begin{center}
    \input{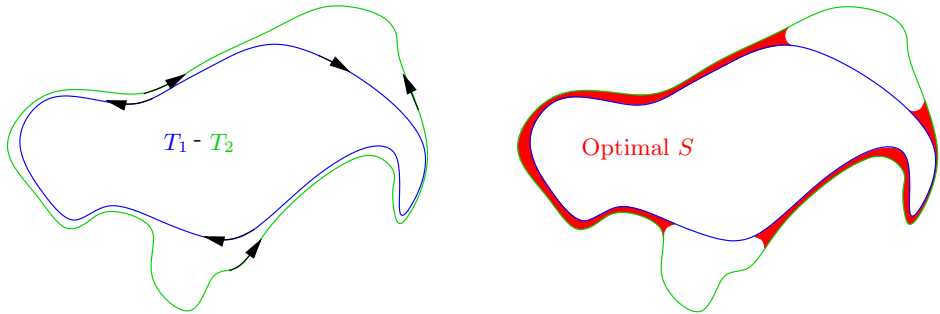}
  \end{center}
  \caption{The flat norm via the $L^1$TV functional provides us with
    both a distance and an informative optimal decomposition into $S$
    and $T_1-T_2 -\partial S$. To use $L^1$TV on sets, we need that
    $T_1 = \partial \Omega_1$ and $T_2 = \partial \Omega_2$ and either
    $\Omega_1 \subset \Omega_2$ or $\Omega_2 \subset \Omega_1$. Using
    a direct method for computing the $\Bbb{F}_{\lambda}$ we don't
    need these inclusions. Alternatively, we can identify the
    difference between two shapes $\partial \Omega_1$ and $\partial
    \Omega_2$ as the boundary of the set
    $\Omega_1\vartriangle\Omega_2$. Then we can use $L^1$TV for all
    codimension 1, boundary shape differences, without requiring
    inclusions.}
  \label{fig:shape-diff}
\end{figure}

As noted above, the flat norm was previously suggested for shape
comparisons in~\cite{glaunes-2005-1, vaillant-2005-1} and then used
in~\cite{glaunes-2007-1} for the purpose of computing shape
statistics. Our observation permits us to use $L^1$TV algorithms to
compute the flat norm distance for many shapes in shape space and the
flat norm decomposition that gives this distance. The decomposition
that we get as a result shows us where the difference is big with
respect to $\lambda$ and where it is small. See
Figure~\ref{fig:shape-diff}.

\section{Summary}
\label{sec:sum}

The innovations introduced in this paper simultaneously expand the
methods available for computing $L^1$TV minimizers, generalizes
$L^1$TV to non-boundary and higher codimensional subsets, opens up 
a new method for multiscale shape decompositions, and supports the
previous suggestion of Glaun\`{e}s et al. that the flat norm is
useful as a distance in shape space.

Difficulties include the fact that non-boundary 0-currents i.e. sets of signed
points in $\Bbb{R}^n$ not arising as the endpoints of a
family of curves, seem rather clumsy to handle. For some applications
the global curvature bound enforced by the method might be too
limiting.  One can imagine a situation in which the curvature of the
approximation should be allowed to vary from place to place on the
input current or set.  One might in fact desire something that returns
a denoised set whose use as a local mean generates a local variance
inversely proportional to the locally allowed curvature.

It is not as easy to handle noise in $T$ in the form of gaps or
missing pieces of $T$. Of course, missing pieces means that $T$ is not
a boundary.  One approach is to first denoise $\partial T$ and then
add the resulting $S$ to $T$ to fill in many or all of the gaps. Then
$T + S$ can be denoised to remove oscillations by computing the
optimal $\Bbb{F}_{\lambda}$ decomposition. But then the denoising
process requires two steps, each involving a choice of $\lambda$ and
this seems a bit clumsy. On the other hand this might make sense
scientifically since the process by which holes and oscillations are
generated may be quite different with different associated length
scales.

In conclusion the program suggested by the relatively simple
observation of the relation between $L^1$TV and the flat norm promises
many new benefits. Many of these benefits are immediately accessible
while others depend on the some further developments outlined above.
These, as well as the general expansion of the above announcement is
the subject of several papers that are in preparation or in planning
with collaborators.
\appendix
\section{Appendix: Micro-tutorial on Currents and the Flat Norm}
\label{sec:mtocfn}

If you know something about currents and have a clear picture of the
flat norm, this section can be skipped. The reference for this section
is Frank Morgan's nice introduction~\cite{morgan-2000-1}. The
definitive, though formidable, treatise for a fair bit of geometric
measure theory is still Federer's 1969 tome~\cite{federer-1969-1}.
References between Morgan and Federer
include~\cite{lin-2002-1,simon-1984-1}.
\begin{description}
\item[Rectifiable sets] Let $\mathcal{H}^n$ denote $n$-dimensional
  Hausdorff measure. An $n$-rectifiable subset $N$ of $\Bbb{R}^{n+k}$
  is the union of 1) an $\mathcal{H}^n$ negligible set and 2) a
  countable collection of subsets of $C^1$ $n$-submanifolds of
  $\Bbb{R}^{n+k}$. We are often interested in the case where
  $\mathcal{H}^n(N) < \infty$.  Intuitively, an $n$-rectifiable set
  looks a great deal like an $n$-manifold at most of its points.
\item[Currents] $n$-Currents in $\mathbb{R}^{n+k}$, denoted
  $\mathcal{D}_{n}(\Bbb{R}^{n+k})$, are the duals to
  $\mathcal{D}^{n}(\Bbb{R}^{n+k})$, the smooth, compactly supported
  $n$-forms on $\mathbb{R}^{n+k}$. We will usually suppress the
  $\Bbb{R}^{n+k}$ and refer simply to $\mathcal{D}_{n}$ and
  $\mathcal{D}^{n}$. We restrict ourselves to integer multiplicity
  rectifiable currents $T$, which have the following representation:
  $T(\phi )=\int_{N}\Theta (x)\left\langle \phi (x),\xi
    (x)\right\rangle d\mathcal{H}^{n},\forall \phi \in
  \mathcal{D}^{n}$ where $N$ is an $n$-rectifiable set in
  $\mathbb{R}^{n+k}$, $\Theta (x)$ is an integer multiplicity density
  function, always $\pm 1$ in this paper, $\phi$ is the form $T$ is
  operating on, and $\xi(x)$ is the unit, simple $n$-vector defining the
  orientation on $N$. Recall that a simple $n$-vector is the wedge
  product of $n$ vectors. In our case, $\xi(x)$ can be thought of as
  an oriented representation of the tangent plane to $N$ at $x$.
  Changing the sign of the density function has the effect of
  reversing the orientation on $N$ which can also be achieved by
  replacing $\xi $ with $-\xi .$

  Currents are naturally equipped with a boundary operator, $\partial
  T(\phi )\equiv T(d\phi )$. $\partial T$ is therefore an
  $(n-1)$-current which operates on $(n-1)$-forms. Note the
  intentional consistency of this definition with Stokes' theorem.

  Intuitively, a current is an oriented manifold or union of oriented
  manifolds, each with an oriented boundary whose orientation is
  inherited from the orientation of the manifold.  See
  Figure~\ref{fig:warped-disk-current}. There are of course wild
  beasts in the menagerie of currents, but unions of $C^1$ manifolds
  with boundary covers a great deal of ground, especially when
  applications are the goal.
\begin{figure}[htp!]
\begin{center}
\input{warped-disk-current.pstex_t}
\end{center}
\caption{Example 2-current $T_M$. Notice that $\partial T_M =
  T_{\partial M}$. The orientation on the boundary $\partial M$ is
  simply that inherited from the orientation on $M$.}
\label{fig:warped-disk-current}
\end{figure}
\item[Mass and the Flat Norm] The \emph{mass} of a current is defined as
  \begin{equation}
    \label{eq:mass-def}
    M(T)=\underset{\left\vert \phi \right\vert \leq 1,\phi \in \mathcal{D}^{n}}{\sup }T(\phi). 
  \end{equation} 
  Informally, the \emph{mass} is simply the $n$-dimensional volume of
  the rectifiable set carrying the $n$-current.  By Theorem 4.1.12 in Federer~\cite{federer-1969-1}, the \emph{flat norm} can be
  defined in two equivalent ways:
  \begin{equation}
    \label{eq:flat-both1}
    \Bbb{F}(T)=\underset{S\in \mathcal{D}_{n+1}}{\min }\left( M(S)+M(T-\partial S)\right) \text{ (given above) }
  \end{equation}
and
\begin{equation}
  \label{eq:flat-both2}
  \Bbb{F}(T) =\underset{\left\vert \phi \right\vert \leq 1,\left\vert d\phi
    \right\vert \leq 1,\phi \in \mathcal{D}^{n}}{\sup }
  (T(\phi ))\text{ (mentioned above) }
\end{equation}
The corresponding  dual definition of the flat norm with scale is given by
\begin{equation}
  \label{eq:flat-scale-dual}
  \Bbb{F}_{\lambda}(T) =\underset{\left\vert \phi \right\vert 
                          \leq 1,\left\vert d\phi
    \right\vert \leq \lambda,\phi \in \mathcal{D}^{n}}{\sup }
  (T(\phi ))
\end{equation}
\item[Examples of the flat norm decomposition] The flat norm involves
  the optimal decomposition of the $n$-current $T$ into an $n$-current
  ($T-\partial S$) and an $(n+1)$-current S. We use the term
  \emph{decomposition} in reference to the fact that $T-\partial S$
  and $S$ are the components explicitly measured by the flat norm,
  even though $T = (T - \partial S) + \partial S$ rather than 
  $T = (T - \partial S) + S$. See Figure~\ref{fig:square-bump-decomp}.
\begin{figure}[htp!]
\begin{center}
\input{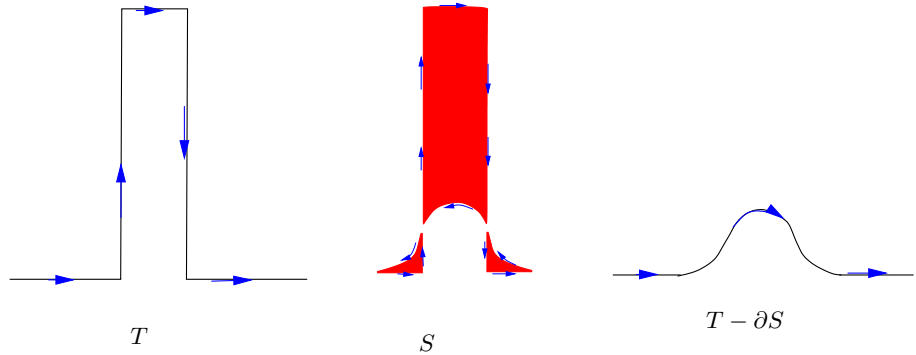}
\end{center}
\caption{Example flat norm decomposition. $T$ is the $1$-current we
  are computing the flat norm of, and $S$ gives the optimal
  decomposition into $S$ and $T-\partial S$. Finally, the flat norm is
  simply $M(S)+ M(T-\partial S)$ = length of the right-most curve and
  area of the red region. }
\label{fig:square-bump-decomp}
\end{figure}
\item[Examples of maximizing forms] The dual formulation of the flat
  norm involves finding the supremum over appropriately constrained
  forms. Figure~\ref{fig:disk-form-field} shows a maximizing form for
  the 2-dimensional disk of radius $r$. We will discuss the
  computation of optimizing forms and the extraction of $S$ from those
  optimizing forms in section~\ref{sec:l1tv-btdf}
  \begin{figure}[htp!]
    \begin{center}
      \input{disk-form-field.pstex_t}
    \end{center}
    \caption{A maximizing form for the disk in 2-D. This form
      satisfies the constraints as long as $\frac{2}{r} \leq \lambda$.
      The $\lambda$ is, of course, the scale in the \emph{flat norm with
      scale} introduced above.}
    \label{fig:disk-form-field}
  \end{figure}
\item[Relation of the flat norm to $L^1$] The flat norm is very much
  like the $L^1$ norm for small differences between currents. The
  value of the flat norm of a difference between two currents ends up
  being roughly the $L^1$ difference between the close parts plus
  the sum of the $n$-volumes of what is left. See
  Figure~\ref{fig:flat-like-l1-plus}.
  \begin{figure}[htp!]
    \begin{center}
      \input{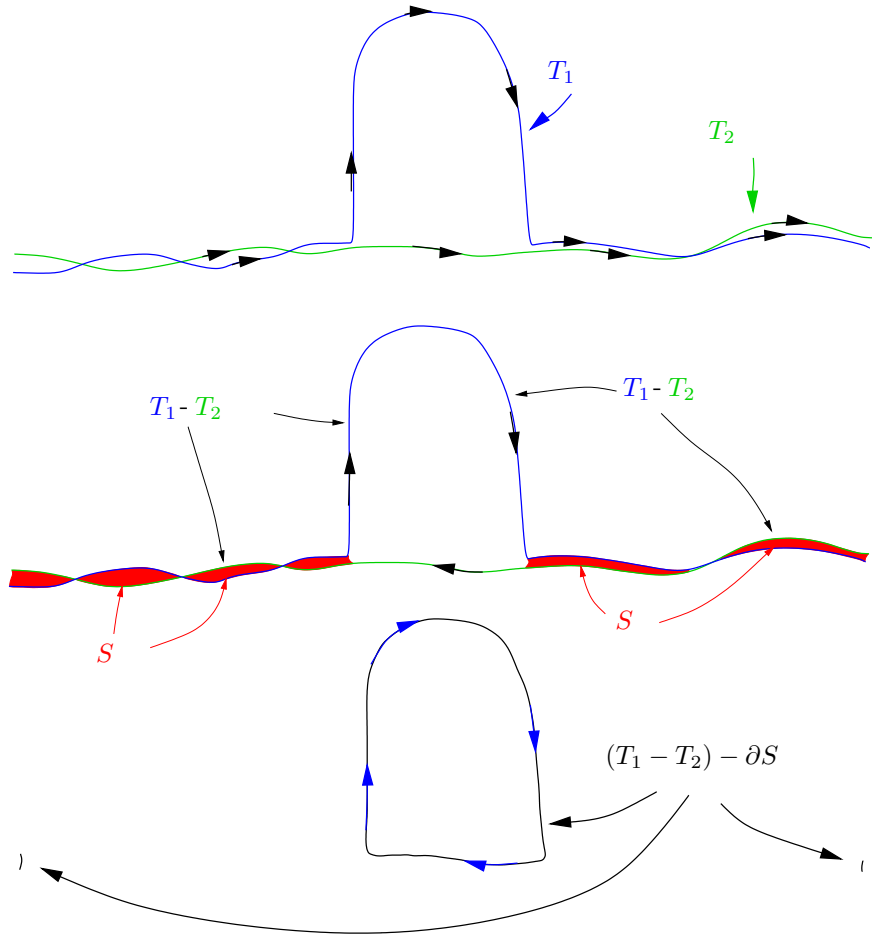}
    \end{center}
    \caption{The flat norm of the difference between these two
      currents is the sum of the area of red region ($L^1$ like) and
      the length of the loop that is left.}
    \label{fig:flat-like-l1-plus}
  \end{figure}
\end{description}

% \bibliographystyle{plain}
% \bibliography{/home/vixie/projects/templates/the_bib}

\end{document}